%Musterdokument: doc.tex
\input amstex
\documentstyle{amsppt}
\magnification=1200
%\input avier %Anweisungen fuer dina4-Format
%Es folgt die "Praeambel", falls erforderlich
\overfullrule=0pt
\NoRunningHeads
%
%\input amacro1
%%%Makros-Dokument%%%allgemeine Makros
%%%%%%%%%%%%%%%%%%%%%%%%%%%%%%%%%%%%%%%%%%%%%
%Blackbold Buchstaben C,F,H,K,N,P,R,Q,Z: z.B. \C
%%%%%%%%%%%%%%%%%%%%%%%%%%%%%%%%%%%%%%%%%%%%%
\def\C{\Bbb C}

\def\N{\Bbb N}

\def\P{\Bbb P}
\def\R{\Bbb R}

%%%%%%%%%%%%%%%%%%%%%%%%%%%%%%%%%%%%%%%%%%%%%%

%
\topmatter
\title
Hinges and geometric constructions of boundaries of Riemann
symmetric spaces
\endtitle
%
%Titel wird fortlaufende Kopfzeile auf den ungeraden Seiten, wenn %nicht
%anders abgegeben
%\rightheadtext{...}
\author
Yu.A. Neretin
\endauthor
%
%
%Autor wird Kopfzeile auf den geraden Seiten, wenn nicht anders angegeben
%\leftheadtext{...}
%\affil 
%Was bedeutet affill?
%\endaffil
%
%
%\address
%Adresse
%\endaddress
\address
Adress Moscow State Institute of electronic and mathematics
\endaddress
%
%
%\curraddr
%derzeitige Adresse
%\endcurraddr
\curraddr
 Max-Planck-Institute fur Mathematik, Bonn
\endcurraddr
%%
%\email
%E-Mail-Adresse
%\endemail
\email
 neretin\@ mpim-bonn.mpg.de ; neretin\@ matan.miemstu.msk.su
\endemail
\date
Datum 9 may 1996
\enddate
\dedicatory
To A.A.Kirilov in his 60 birthday
%Widmung
\enddedicatory
\thanks
I would like to thank C.De Concini, S.L.Tregub, M.M.Kapranov, 
W.Ballmann and M.A.Olshanetsky for discussion of this subject.
I also grateful to administrations of Max-Plank-Institute
fur Mathematik (Bonn), Mittag-Leffler Institute (Stokgolm)
and Erwin Shrodinger Institute (Wien) for hospitality.
\endthanks
%
%
%\translator
%Uebersetzer
%\endtranslator
%
%
%
%
%\subjclass
%Klassifizierung
%\endsubjclass
%
%
%
\abstract
 We give elementary constructions for Satake-Furstenberg, Martin and
Karpelevich  boundaries of symmetric spaces.  We also consruct some 
"new" boundaries
\endabstract
%\input toc...
%Muster fuer Inhaltsverzeichnis: toc.tex
%\toc
%oder
%\toc\nofrills{Eigene Ueberschrift}
%\widestnumber\specialhead{...}
%\widestnumber\head{...}
%\widestnumber\subhead{...}
%\widestnumber\subsubhead{...}
%\head {..} ... \page{..}
%\endhead
%
%\endtoc
\endtopmatter
\document
%\hyphenation{...}

It is well-known that symmetric spaces have nontrivial and
nice  boundaries. There are two
(disjoint) scientific traditions of investigation of such boundaries.
The first tradition is related to enumerative
algebraic geometry of quadrics. It was begun by the paper of Study
(1886) on the geometry of the space $PGL(3,\C)/SO(3,\C)$
(this is the space of all nondegenerate conics in $\C\P^2)$. This
construction was extended by Semple (1948-1951) to the spaces
$PGL(n,\C)/SO(n,\C)$ and to the groups $PGL(n,\C)$ itself. Later 
(1983)
De Concini and Procesi constructed analogical compactification for
arbitrary symmetric space $G/K$ where  $G$ is a semisimple group 
without center and $K$ is a {\it complex} symmetric subgroup.

Another scientific tradition is related to harmonic analysis on symmetric 
spaces. In 1960 Satake constructed nice compactifications of Riemanian 
symmetric spaces( these compactifications
 are real forms of compactifications of complex 
symmetric spaces mentioned above). In 1961-1969 in  Karpelevich, Dynkin
and Olshanetsky constructed more complicated
boundaries(their works were devoted to analysis of harmonic functions
on the symmetric spaces).

The purpose   of these notes (it is a part of the paper [33])
 is to give elementary description
for Satake-Furstenberg boundary, Karpelevich boundary,
Martin (Dynkin-Olshanetsky) boundary for Riemann noncompact symmetric spaces,
we also construct some "new" boundaries(velocity boundaries
in section 3 and sea urchins in section 6).
We discuss only the boundaries of symmetric spaces $PGL(n,\R)/SO(n,\R)$
(boundaries of other classical symmetric spaces can be described by 
the same way).

%I thanks C. De Concini, S.L.Tregub, M.M.Kapranov, W.Ballmann and 
%M.A.Olshanetsky for discussion of this subject.

\head{1.Preliminaries. Hinges}
\endhead

\demo{1.1. Linear relations }
Let $V,W $ be  linear spaces. A {\it linear relation} $V\rightrightarrows W$ 
is arbitrary linear subspace in $V\oplus W$.

\demo{Example 1.1.} Let $A:V\rightarrow W$ be a linear operator.
Then its graph $graph(A)$ is a linear relation.
\enddemo

Let $P:V\rightrightarrows W$   be a linear relation. Then we define

1. the {\it kernel} $Ker(P)=P\cap (V\oplus 0)$

2. the {\it image} $Im(P)$ is the projection of $P$ to $0\oplus W$

3. the {\it domain} $Dom(P)$ is the projection of $P$ to $V\oplus 0$

4. the {\it indefinitness} $Indef(P ) = P\cap (0\oplus W)$

\demo{Remark 1.2} Let $P=graph(A)$. Then $Im(P)$ is the usual
image of the linear operator $A$ and $Ker(P)$ is the usual kernel
of the linear operator $A$.

\enddemo

We also define the {\it rank} of a linear relation $P$:
$$
rk(P)=\dim Dom(P)-\dim Ker(P)=\dim Im(P)-\dim Indef(P)= $$  $$
= \dim P\ -
\dim Ker(P)-\dim Indef(P)
$$

 \demo{Remark 1.3} Let us consider a linear relation $P:V\rightrightarrows W$.
Then it defines by the obvious way the invertible linear operator
$$
[P]:Dom(P)/Ker(P)\rightarrow Im(P)/Indef(P)
$$

 \enddemo

\demo{1.2. Nonseparated quotient of grassmanian}
 We denote by $\R^\ast$ the multiplicative group of $\R$.
We denote by $Gr_n$ the grassmanian of all $n$-dimensional subspaces in 
$\R^n\oplus \R^n$. Let $P:\R^n\rightrightarrows R^n$ be a element
of $Gr_n$. Let $\lambda\in\R^\ast$. We define $\lambda \cdot P\in Gr_n$
by the condition
$$
(v,w)\in P \Leftrightarrow (v,\lambda w)\in \lambda P
$$
\demo {Remark 1.4} If $P$ has the form $graph(A)$
then $\lambda\cdot P=graph(\lambda A)$.
\enddemo

 Let us consider the quotient space $Gr_n/\R^\ast$ equipped with the usual
quotient topology (see [29]).
Let us consider a sequence $x_j\in Gr_n/\R^\ast$ and a point
$y\in Gr_n/\R^\ast$. Let $P_j,Q$ be representatives of $x_j$ and $y$
in $Gr_n$. Then the sequence $x_j$ converges to $y $ if there exist
$\lambda_j,\lambda\in \R^\ast$ such that $\lambda_j\cdot P_j$ 
converges to $\lambda\cdot Q$ in the topology of $Gr_n$.

We will use the same notations for points $P \in Gr_n$ and their 
$\R^\ast$-orbits, i.e. we denote the orbit $\R^\ast\cdot P$ by $P$.

 There are two types of orbits of $\R^\ast$ on
$Gr_n$. If $rk(P)=0$(i.e $P=Ker(P)\oplus Indef(P)$) then $P$ is a 
fixed point of the group $\R^\ast$. If $rk(P)\not=0$ then the 
stabilizer of $P$ in $\R^\ast$ is trivial and hence the orbit is 
isomorphic to the group $\R^\ast$ itself. The orbit of the first type 
are closed. The orbits of   the second type are not closed. Hence 
the topology in the space $Gr_n/\R^\ast$ is not separated in the 
Hausdorff sence. A point $P\in Gr_n$ is closed set only in the case
$rk(P)=0$.

\demo{Example 1.5} Let us consider a sequence
$$
A_j=\pmatrix j&0\\ 0&1\endpmatrix
$$
of linear operators in $\R^2$. Let  
$P^j\in Gr_2$ be their graphs. Let us consider the 
sequences
$$
j\cdot P_j\ ;\ P_j \ ;\ j^{-1/2}\cdot P_j\ ;\ j^{-1}\cdot P_j\ ;\ 
j^{-2}\cdot P_j
$$
in $Gr_2$. Their limits in $Gr_2$ are the subspaces  $R_1,\dots,R_5$ 
having the form
$$
R_1:(0,0;x,y)
$$
$$
R_2:(0,y;x,y)
$$
$$
R_3:(0,y;x,0)
$$
$$
   R_4:(x,y;x, 0)
$$
$$
R_5:(x,y;0,0)
$$
Hence the sequence $P_j$ has 5 limits in the quotient space
$Gr_2/\R^\ast$.
\enddemo
\demo{Remark 1.6}
Let us consider a sequence of invertible operators 
$A_j:\R^n\rightarrow \R^n $. 
Let $P_j$ be their graphs.
Evidently subspaces $\R^n\oplus 0$ and $0\oplus \R^n$ are limits
of the sequense $P_j$ in the quotient space
$Gr_n/\R^\ast$            .
By the official topological definition this sequence is convergent
(and moreover it has at least 2 limits). It is quite clear
that official definition of convergence ( the sequence converges if 
it has limit) is bad.
\enddemo

 Let $A_j$ be a sequence of invertible operators. Let $P_j$ be their
graphs. We say that the sequence $P_j$ is {\it seriously convergent} if
each limit point  of $P_j$ in  the quotient space
$Gr_n/\R^\ast$ is the limit the limit 
          of $P_j$ in  the quotient space.

\demo{Remark 1.7}
We define serious converegrnce only for sequences of {\it invertible}
operators!
\enddemo

We say that the subset $S\in Gr_n/\R^\ast$ is {\it admissible} if
there exists seriously convergent sequence $P_j$ such that the
set of limits of $P_j$ coincides with $S$.

\demo{Example 1.8} The sequence $P_j$ described in example 1.5  is 
seriously convergent. Hence the set $R_1,\dots R_5$ is admissible.
\enddemo

\enddemo

\demo{1.3. Hinges}

\demo{Definition 1.9}
 A {\it hinge} $$\goth P=(P_1,\dots,P_k) $$
is a family of elements of $Gr_n/\R^\ast$ such that

$0^\circ$. For all j $rk(P)>0$

$1^\circ$. For all $j$
 $$Ker(P_j)=Dom(P_{j+1})$$
$$Im(P_j)=Indef(P_{j+1})
                    $$

$2^\circ$. $$Indef(P_1)=0$$ 
   $$Ker(P_k)=0$$
i.e. $P_1$ is the graph of a operator 
$(\R^n\oplus 0)\rightarrow (0\oplus\R^n)$
and $P_k$ is the graph of a operator $(\R^n\oplus 0)\leftarrow (0\oplus\R^n)$

We denote space of all hinges in      $\R^n$ by $Hinge(n)$
\enddemo
\demo{Remark 1.10}     The condition $2^\circ$ is intepretation of 
the
condition $1^\circ$ if $j=0$ and $j=k$.
\enddemo
\demo {Example 1.11} The graph of a invertible operator is a hinge 
(k=1). The graph of a noninvertible operator is not hinge (see the 
condition $2^\circ$)
                 \enddemo
\demo{Example 1.12} In  Example 1.4 the set  
       $$(R_4,R_2)$$
is a hinge. Note that the rank of $R_1,R_3,R_5$ is $0$.
                     \enddemo

By the definition af hinge we have
$$
\R^n= Dom(P_1)\supset Ker(P_1)=Dom(P_2)\supset Ker(P_2)=\dots
$$
Hence (by the condition $0^\circ$)  we have $k\le n-1$

\demo{Theorem 1.13} {\it Let us consider a hinge
$$\goth P=(P_1,\dots,P_k)$$
Let
$$Q_j=Ker (P_j)\oplus Im(P_j)= Dom(P_{j+1})\oplus Indef 
(P_{j+1})\in Gr_n/\R^\ast$$
$$Q_0=\R^n\oplus 0$$
Then the set
$$\{Q_0,P_1,Q_1,P_2,\dots,P_k,Q_k\}$$
is a admissible subset in $Gr_n/\R^\ast$ .  Moreover each
admissible subset has such form.}

\enddemo
\demo{Remark 1.14}
Unformally speaking hinges are limits in $Gr_n/\R^\ast$ 
of sequences of 
invertible operators. For 
instance sequence $A_j$ described in the Example 1.  converges to the
hinge $(R_4,R_2)$ .
Hinges are slightly different from admissible sets. Neverless it is
better for us to forget about fixed points $Q_0,Q_1,\dots$ (since
they can be reconstructed by $P_1,P_2,\dots$)
\enddemo
\enddemo
 \demo{1.4. The topology on the space of hinges} 
  Let $M$ be a compact metric space with a metric
$\rho(\cdot,\cdot)$. Let $\Cal S(M)$ be the space of 
  all closed subsets in $M$. Let $X\in \Cal S(M)$
. We denote by $\Cal O_{\epsilon}(X)$ the set of
points $m\in M$ such that exists $x\in X$ satisfying the condition
$\rho (m,x),\epsilon$.

Let $X_1,X_2 $ be closed subsets.
{\it Hausdorff distance} (see [31]) 
between $X_1$ and $X_2$  is infimum of $\epsilon$
such that
$$
X_1\subset \Cal O_\epsilon(X_2)\ \ ;\ \ X_2\subset \Cal O_\epsilon(X_1)
$$
It is well known that the space $\Cal S(M)$ equipped with the Hausrorff
distance is a compact metric space.

Let us consider a invertible operator $A$ and its graph $P$. Let us 
consider the curve $\R^\ast\cdot P$ in grassmanian.
Let us consider its closure $\sigma(A)$. It contains the curve
$\R^\ast\cdot P$ itself and two points $\R^n\oplus 0$ ,
$0\oplus\R^n$.        We denote family of curves 
$\sigma(A)\in \Cal S(Gr_n)$ by $\widetilde{PGL(n,\R)}$. We have 
the obvious bijection
$$
{PGL(n,\R)}\leftrightarrow \widetilde{PGL(n,\R)}
$$
 We denote by $\overline{PGL(n,\R)}$
   the closure of $\widetilde{PGL(n,\R)}$
in the Hausdorff metric.

\demo{Theorem 1.15}{\it 
Let
$$\goth P=(P_1,\dots,P_k)$$
be a hinge. Let $Q_j$ be the same as in the theorem 1.13. Let us denote by
$\gamma(\goth P)$ the curve
                           $$
Q_0\cup(\R^\ast\cdot P_1)\cup Q_1\cup(\R^\ast\cdot P_2)\cup Q_2\dots\cup 
Q_k
                               $$
Then the map
$$
 \goth P\mapsto \gamma(\goth P)
$$
 is the bijection
                 $$
 Hinge(n)\rightarrow              \overline{PGL(n,\R)}
$$             }       
                       \enddemo

We see that $Hinge(n)$ has the natural structure of a
 compact metric (metrizable) space containing $PGL(n,\R)$
as open dence set( if $A\in PGL(n,\R)$ then its graph 
is a one-element hinge $\goth P =(P)$
).

 The space $Hinge(n)$ has the natural
structure of $(n^2-1)$-dimensional real analytic manifold (it is not 
obviuous).
The set $Hinge(n)\backslash PGL(n,\R)$ is the union of $(n-1)$ 
submanifolds  of codimension 1 (see below bibliographical remarks).

\enddemo

\enddemo

\head {2. Satake-Furstenberg boundary}
        \endhead
\demo{2.1.Symmetric space $SL(n,\R)/SO(n)$}
Let us consider the 
space $Q$ of  real symmetric  positive definite matrices defined
up to multiplier. The action of the group $SL(n,\R)$ on this
space is defined by the formula
$$
 g:A\mapsto g A g^t
$$
where $A$ is symmetric matrix, $g\in SL(n,\R)$ and $g^t$ is the 
transposed matrix. Obviously the stabilizer of the point $E$ is
the group $SO(n)$. Hence we obtain
$$
  Q\simeq SL(n,\R)/SO(n)
                      $$
\enddemo
        \demo{2.2. Positive linear relations}We want to describe
the closure of the space $Q$ in $Hinge(n)$. For this purpose
we need in some 
preliminaries.
  Let us consider in the space $\R^n$ the standard scalar product
$$
<v,w>=\sum_k v_kw_k
$$
We define in the space $\R^n\oplus \R^n$ the skew-symmetric bilinear form
by the formula
$$
\{(v,v');(w,w')\}       = <v,w'>-<w,v'>
$$
We define also indefinte symmetric bilinear form on  $\R^n\oplus \R^n$
     by the 
formula
$$
[(v,w);(v',w')]:=<v,w'>+<v',w>
$$
We say that a $n$-dimensional linear relation 
$P:\R^n\rightrightarrows\R^n$ is {\it symmetric} if
$P$
  is a maximal isotropic subspace with respect to the skew-symmetric
bilinear form $<\cdot,\cdot>$.
 
 \demo{Remark 2.1} Let $A  $ be a symmetric linear operator(i.e 
 $A=A^t$). Then its graph is a symmetric linear relation.
\enddemo

  Let us consider a symmetric 
linear relation $:\R^n\rightrightarrows \R^n$  
. Then $Im(P)$ is the orthogonal complement in $\R^n$ to $Ker(P)$
(with respect to the scalar product $<\cdot,\cdot>$) 
and $Indef(P)$ is the orthogonal complement to $Dom(P)$ (with respect to
the standard scalar product in $\R^n$ ). Hence the linear relation $P$
 defines the nondegenerate pairing
$$
Dom(P)/Ker(P)\times Im(P)/Indef(P)\rightarrow \R                 \tag 2.1
$$
The linear relation $P $ also defines the operator
$$
Dom(P)/Ker(P)\rightarrow Im(P)/Indef(P)                 \tag 2.2
$$
Hence each symmetric linear relation $P$ defines nondegenerated 
symmetric bilinear form $\goth q_P$ on the space $Dom(P)/Ker(P)$.

We say that a symmetric linear relation $P$ is {\it nonnegative 
definite} if 
 the form
$[\cdot,\cdot]$ is nonnegative definite on the 
subspase $P$.  It is equivalent 
to 
the positivity of quadratic form $\goth q_P$.
\demo {Remark 2.2} Let a linear relation $P$ be the graph of a 
operator A. Then $P $ is nonnegative definite if and only if $A$ is
nonnegative definite.
\enddemo
\enddemo
\demo{2.3.Satake-Furstenberg boundary}      Let us consider the closure
$\overline Q$ of the space $Q$ in the space $Hinge( n)$.
It is easy to show that a hinge $\goth P$ belongs to $\overline Q$
if and only if all linear relations $P$ are nonnegative definite.
It appears that 
this closure coincides with 
  {\it Satake-Futrsenberg compactification}  of the symmetric space
$SL(n,\R)/SO(n)$.

Hence a point of Satake-Furstenberg compactification 
is given by the following  data:

  $1^\ast$.$ s=1,2,...,n-1$

  $2^\ast$. A hinge
$$\goth P=(P_1,\dots,P_s)$$ such that    
all linear relations $P_j$ are nonnegative definite.

  Let us consider a point of Satake-Furstenberg compactification (i.e data 
  $1^\ast-2^\ast$) . Let us consider the subspaces
$$
V_j=Ker(P_j)=Dom(P_{j+1})
$$
Then the form $\goth q_{P_j}$ is positive definite 
on      $Dom(P_j)/Ker(P_j)$ . Now we can say that a point of 
Satake-Furstenberg boundary is defined by the following  data

  $1^\star$.$ s=1,2,...,n-1$

  $2^\star $. A flag
$$
0\subset V_1\subset V_2\subset\dots\subset V_s  \subset \R^n
$$
where all subspaces $0,V_1,\dots, V_s, \R^n$ are different  .

$3^\star$.A positive definite quadratic form $R_j$ in each quotient 
space $Dom(P_j)/Ker(P_j)$.

\enddemo

\head{3. Velocity compactifications of symmetric spaces.}\endhead

  \demo{3.1. Simplest velocity compactification}\enddemo
 Let $A\in Q=SL(n,\R)/SO(n)$
be a positive definite matrix. Let
$$a_1 \ge  a_2 \ge ...\ge  a_n$$
be eigenvalues of $A$. Let
 $$
  \lambda_j = \ln a_j
$$
We denote by $\Lambda(A)$ the collection
$$ \align \Lambda (A) = (\lambda_1, \lambda_2,
 ... ,\lambda_n)\ ; \ \ \lambda_1 \ge
\lambda_2 \ge ... \tag 3.1\endalign $$
 The matrix $A$ is defined up to multiplier
  and hence $\Lambda(A)$ is defined
up to additive constant:
$$\align 
(\lambda_1, \lambda_2, ... \lambda_n) 
\sim (\lambda_1 +\sigma, \lambda_2+\sigma, ... ,\lambda_n+\sigma)
\tag 3.2 \endalign$$
We denote by $\Sigma_n$ the space 
of all collections $\Lambda (A)$(see (3.2)).
It is easy to see that $\Lambda(A)$ is a $(n-1)$-dimensional siplicial cone.
We can assume $\lambda_n=0$ and hence the cone $\Sigma_n$ 
can be considered as the space
of collections
$$
\{\lambda_1 \ge \lambda_2 \ge ...\ge \lambda_{n-1}\ge 0\} 
$$
We denote by $\Delta_n = \partial \Sigma_n $ the $(n-2)$-dimensional simplex 
$$
1 \ge \mu_2 \ge \mu_3 \ge ...\ge \mu_{n-1} \ge 0
$$
  It is natural to  think that $\mu_1=1,\mu_n=0$.
We say that $\Delta_n$ is the {\it velocity simplex}.
Let us consider the natural projection
$$
\pi : (\Sigma_n \backslash 0 )\rightarrow \Delta_n
$$
defined by the rule
$$
\pi(\lambda_1, \lambda_2, ... \lambda_{n-1}, 0)=
({\lambda_2 \over \lambda_1}, {\lambda_3 \over \lambda_1},
..., {\lambda_{n-1}\over
\lambda_1})
$$
Now we define the compactification
$$\overline \Sigma_n=\Sigma_n \cup \Delta_n$$
of $\Sigma_n$.
A sequence $L_j =(\lambda_1^{(j)},...,\lambda_n^{(j)})\in \Sigma_n$ 
converges
to $M\in \Delta_n$ if

  1. $\lambda_1^{(j)} - \lambda_n^{(j)} \rightarrow
 \infty$ if $j \rightarrow \infty$

  2.The sequence $ \pi(L_j) \in \Delta_n$ converges to $M$.

  We also define the {\it velocity compactification}
$$
\overline Q^{vel} =SL(n,\R)/ SO(n) \cup \Delta_n$$
of the symmetric space $SL(n,\R)/ SO(n)$. 
A sequence $A_j$ in $Q$ converges to $M\in \Delta_n$ if
 $\Lambda(A_j)$ converges to $M$ in the topology of $\overline\Sigma_n$.

 \demo{3.2. Polyhedron of Karpelevich velocities}\enddemo
Now we want to describe more delicate compactification of the simplicial 
cone
$\Sigma_n$ (compactification by Karpelevich velocities). Let us consider a 
sequence
$$\lambda^{(j)}=\{\lambda_1^{(j)}\ge\dots\ge\lambda_n^{(j)}\}\in\Sigma_n$$
Let $1\ge\mu_2\ge\dots\ge\mu_{n-1}\ge 0$ be its limit in $\Delta_n$.
It can happens that some of numbers $\mu_i$ are equals:
$$
\mu_k=\mu_{k+1}=\dots=\mu_l
$$
In this case we will separate velocities of
$$
 \{\lambda_k^{(j)}\ge\dots\ge\lambda_l^{(j)}\} \in \Sigma_{l-k+1}
                                      $$
by the same rule as above. 
%Formal definition is complicated. 
%Nonformally the definition is quite simple and it is illustrated by
%the Example 3.8   below.

%%%%%%%%%%%%%%%%%%%%%%%%%%%%%%%%

\demo{Definition of the polyhedron}
We denote by $I_{\alpha,\beta}$   set $\{ 
\alpha,\alpha+1,\dots,\beta\}\subset \N$

%%%%%%%%%%%%%%%%%%%%%%%%%%%%%%%%%%%%%%
Let us consider a interval 
$I_{\alpha,\beta}= \{\alpha, \alpha+1,\dots,\beta\}$
. We denote by $\Sigma(I_{\alpha,\beta})$ the simplicial cone
$$
\lambda_\alpha\ge\lambda_{\alpha+1}\ge\dots\ge\lambda_\beta
$$
the elements of the cone $\Sigma(I_{\alpha,\beta})$
  are defined up to additive constant (see (3.2)).
 We also define  the simplex $\Delta(I_{\alpha,\beta})$ given by the 
 unequalities
$$
1=\mu_{\alpha}\ge\mu_{\alpha+1}\ge\dots\ge\mu_{\beta-1}\ge\mu_\beta=0
$$
Let us consider the compactification
$$
\overline\Sigma(I_{\alpha,\beta})=\Sigma (I_{\alpha,\beta})\cup
\Delta(I_{\alpha,\beta})
$$

\demo{Remark 3.1 } Let us consider the case 
$\alpha=\beta$. The set 
$\Sigma(I_{\alpha,\alpha})=\overline\Sigma(I_{\alpha,\alpha})$ 
consist of the unique point (it is one 
real
number defined up to additive constant).

%%%%%%%%%%%%%%%%%%%%%%%%%konets kuska

%%%%%%%%%%%%%%%%%%%%%%%%nachalo vstavki
Let $k\le \alpha\le\beta\le l$. We define the map
$$
\Pi^{k,l}_{\alpha,\beta}:\Sigma(I_{k,l})\rightarrow \Sigma(I_{\alpha,\beta})
$$
given by the formula
                        $$
\Pi^{k,l}_{\alpha,\beta}(\lambda_k,\dots,\lambda_l)=(\lambda_\alpha,\dots
,\lambda_\beta)
$$
We define two polyhedra
$$
\Xi(k,l):= \prod_{\alpha,\beta:k\le\alpha\le l\le\beta}
\Sigma(I_{\alpha,\beta})
$$
$$
\overline\Xi(k,l):= \prod_{\alpha,\beta:k\le\alpha\le l\le\beta}
\overline\Sigma(I_{\alpha,\beta})
$$
Obviously         $\Xi(k,l)\subset\overline\Xi(k,l)$.
Let us consider the natural (diagonal )embedding
$$
i:\Xi(I_{k,l})\rightarrow             \Xi(k,l)
$$
(it is the product of the maps  $\Pi^{k,l}_{\alpha,\beta}$)

  The {\it polyhedron of Karpelevich velocities}  $\Cal(k,l)$ 
is the closure of
the set $i(\Sigma(I_{k,l}))$ in      $\overline\Xi(k,l)$.
\enddemo

\demo{Criterium of convergence of a sequence of 
interior points to a point of the 
boundary}

Let us consider a sequence
$$
 \Lambda^{(j)}=\{\lambda^{(j)}_k,\lambda^{(j)}_{k+1},\dots,\lambda^{(j)}_l\}
$$
 Then the nessesary and sufficient condition of convergence of
the sequence $\Lambda^{(j)}$   in $\Cal K(k,l)$ is the convrgence of
all sequences
$$
\Pi^{k,l}_{\alpha,\beta}( \Lambda^{(j)})
=(\lambda^{(j)}_\alpha,\dots,\lambda^{(j)}_\beta)
$$
in $\overline\Sigma(I_{\alpha,\beta})$.

The Karpelevich velocity polyhedron is defined. Now we want to
give explicit description of its combinatorical structure.
 \enddemo

%%%%%%%

%
%%%%%%%%%%%%%%%%%%%%%%%%%%konets vstavki

\demo{Tree-partitions}
   Let us consider the set $I_{k,l}:=\{k,k+1,...,l\}$.
We say that a {\it partition} of $I_{k,l}$ is a representation of $I_{k,l}$
as
$$
I_{k,m_{1}}\cup I_{m_1+1,m_2}\cup\dots\cup I_{m_{s-1}+1,l}
$$
where $s>1$.

  We say that a system
$\goth a$ of subsets of $I_{k,l}$ is a {\it tree-partition} if

  a) $I_{k,l}\in \goth a$ 
   
  b) Each element $J \in \goth a$ has the form
 $I_{\alpha, \beta}=\{\alpha,\alpha+1,\dots,\beta\}$

   c) If $J_1,J_2\in \goth a$ then

$$J_1\cap J_2=\emptyset \ \ or \ \ J_1\supset J_2 \ \ or \ \ J_2\subset J_1$$

   d) Let $J=I_{\alpha, \beta}\in \goth a$ .Then there are
 only two possibilities

   \ \ \ \ $1^\ast$.There is no $K\in \goth a $ such that $K\subset
 J$ (in this case we say that $I_{\alpha, \beta}$ is  
 {\it irreducible}).

\ \ \ \ $2^\ast$. $J=I_{\alpha, \beta}$ can be decomposed as the union

$$
I_{\alpha, \beta}=I_{\alpha, \gamma_1}\cup I_{\gamma_1 +1,\gamma_2}\cup
I_{\gamma_2 +1,\gamma_3}\cup ...\cup I_{\gamma_{s-1} +1,\beta} \tag 3.3$$
where $I_{\alpha, \gamma_1}, I_{\gamma_1 +1,\gamma_2},...,
 I_{\gamma_{s-1} +1,\beta}\in \goth {a}  $.
In this case we say that $J$ is {\it reducible} and (3.3) is the 
{\it canonical decomposition} of $J$.
  \demo{Remark 3.2}Let $I_{\alpha, \beta}\in \goth a$.
 Let $\goth b$ be the
set of all $J\subset I_{\alpha, \beta}$ such that $J \in \goth a$. 
Then $\goth b$ is the tree-partion of $I_{\alpha, \beta}$.\enddemo

  \demo{Remark 3.3} In the other words  tree-partition
 is given by the following  data.
 We consider a partition of the segment
 $I_{k,l}\subset \Bbb N$ to subsegments, then
 we consider partitions of some subsegments,  etc.\enddemo

We denote by $TP(k,l)$ the set of all tree-partitions of $I_{k,l}$. Let 
us define the partial canonical ordering on
$TP(k,l)$. Let $\goth a,\goth b\in TP(k,l)$. We say that $\goth a      
        >\goth b$ if $J\in \goth a$ implies $J\in \goth b$ (i.e
$\goth b \supset \goth a$).

The partially ordered set $TP(k,l)$ contains the unique
maximal element $\goth a_0$. This is the tree-partition which 
contains the unique element $I_{k,l}$.

  A element $\goth b\in TP(k,l)$ is minimal if

a) Each irreducible element of $\goth b$ contains only one point.

b)If $J\in \goth b$ is reducible then the canonical decomposition of
$J$ contains   exactly two elememts ($s=2$ in (3.3)  ).

\enddemo

\demo{Description of the polyhedron}
%Let us consider a interval 
%$I_{\alpha,\beta}= \{\alpha, \alpha+1,\dots,\beta\}$
%. We denote by $\Sigma(I_{\alpha,\beta})$ the simplicial cone
%$$
%\lambda_\alpha\ge\lambda_{\alpha+1}\ge\dots\ge\lambda_\beta
%$$
%the elements of the cone $\Sigma(I_{\alpha,\beta})$
%  are defined up to additive constant (see (3.2)).
% We also define  the simplex $\Delta(I_{\alpha,\beta})$ given by the 
% unequalities
%$$
%1=\mu_{\alpha}\ge\mu_{\alpha+1}\ge\dots\ge\mu_{\beta-1}\ge\mu_\beta=0
%$$
%Let us consider the compactification
%$$
%\overline\Sigma(I_{\alpha,\beta})=\Sigma (I_{\alpha,\beta})\cup
%\Delta(I_{\alpha,\beta})
%$$
%
%\demo{Remark 3.3 } Let us consider the case 
%$\alpha=\beta$. The set 
%$\Sigma(I_{\alpha,\alpha})=\overline\Sigma(I_{\alpha,\alpha})$ 
%consist of the unique point (it is one 
%real
%number defined up to additive constant).
%
%\enddemo

Let us consider a partition $\goth t$ of $I_{\alpha,\beta}$:
$$
I_{\alpha, \beta}=I_{\alpha, \gamma_1}\cup I_{\gamma_1 +1,\gamma_2}\cup
I_{\gamma_2 +1,\gamma_3}\cup ...\cup I_{\gamma_{s-1} +1,\beta} \tag 3.4$$
We denote by $\widetilde\Delta(I_{\alpha,\beta}|\goth t)$ the open simplex
$$
{1=\mu_\alpha=\dots=\mu_{\gamma_1}>\mu_{\gamma_1+1}=\mu_{\gamma_1+2}=\dots
=\mu_{\gamma_2}>\dots>
\mu_{\gamma_{s-1}+1}=\dots=\mu_\beta=0}          \tag 3.5
$$
We denote by $\Delta(I_{\alpha,\beta}|\goth t)$ the compact simplex
$$
{1=\mu_\alpha=\dots=\mu_{\gamma_1}\ge\mu_{\gamma_1+1}=\mu_{\gamma_1+2}=\dots
=\mu_{\gamma_2}\ge\dots\ge
\mu_{\gamma_{s-1}+1}=\dots=\mu_\beta=0}          \tag 3.6
$$
It is natural to consider in $\Delta(I_{\alpha,\beta}|\goth t)$ and 
$\widetilde\Delta(I_{\alpha,\beta}|\goth t)$ the coordinates
$$       \tau_2:=\mu_{\gamma_1+2}=\dots
=\mu_{\gamma_2}=
$$
$$
\dots\dots
$$
$$
\tau_{s-1}:=\mu_{\gamma_{s-2}+1}=\dots=\mu_{\gamma_{s-1}}
$$
\demo{Remark 3.4} If s=2 then $\Delta (J\ |\goth t)=
\widetilde\Delta (J\ |\goth t)$ consist of the unique
point $\{1>0\}$.
\enddemo

\demo{Remark 3.5}
$$
\Delta(I_{\alpha,\beta})=\bigcup_{\goth t}\widetilde\Delta (I_{\alpha,\beta}\ 
|\goth t)
$$
where the union is given by the all partitions of $I_{\alpha,\beta}$
\enddemo

Fix a tree-partition $\goth a\in TP(k,l)$. For each element 
$J\in\goth a$ consider its canonical decomposition $\goth t$. We 
denote the
simplex $\widetilde\Delta(J\ |\goth t)$  by 
$\widetilde\Delta(\goth a,J)$ .

For each          $\goth a\in TP(k,l)$ we define the {\it face}
$$
F(\goth a)=
(\prod_{J=I_{\alpha,\beta}\in\goth a \ is \ irreducible}
\Sigma(I_{\alpha,\beta}))\times
$$
$$
\times\prod_{J\in\goth a \ is \ redicible}
\widetilde\Delta(\goth a, J)                     \tag 3.7
$$
 \demo{Remark 3.6} For the trivial tree-partition $\goth a_0$ we
have $F(\goth a_0)=\Sigma(I_{k,l})$. If $\goth b $ is a minimal tree 
-partition then $F(\goth b)$ is a one-point-set.
\enddemo

  We define {\it Karpelevich velocity polyhedron} $\Cal K({k,l})$ 
  by
$$
\Cal K({k,l})=\bigcup_{\goth a\in TP(k,l)} F(\goth a)
$$
We want to define a topology of a compact metric space on 
$\Cal K({k,l})$. The face $F(\goth a_0)=\Sigma(I_{k,l})$ will be a 
open dense subset 
in             $\Cal K({k,l})$.
\demo{Remark 3.7}
Let $l=k$. Then $\Cal K(k,k)$ consist of one point. Let $l=k+1$. Then 
we have two tree-partitions of the set $\{k,k+1\}$: The trivial
tree-partition $\goth a_0$ and maximal tree-partition $\goth a_1$
(its elements are $(k,k+1),(k),(k+1)$).The face $F(\goth a_0)$ is closed
half-line $\lambda_1>0$                . The face
$F(\goth a_1)$ is one-point-set. Hence $\Cal K (k,k+1)$ is the 
segment $[0,\infty]$.
\enddemo
\enddemo

\demo{Convergence of interior points to the boundary}

The definition of convergence is inductive. We assume  the 
convergence is defined for all Karpelevich polyhedra 
$\Cal K({\alpha,\beta})$ such that $\beta-\alpha<l-k$.

    We define the convergence of a  sequence
$$
x^{(j)}=\{x_k^{(j)}\ge\dots\ge x_l^{(j)}\}\in\Sigma(I_{k,l})=F(\goth a_0)
$$
in two steps.

{\it The first step.} The convergence of $x^{(j)}$ in $\overline\Sigma(I_{k,l})$
is a nessesary condition for the convergence in    $\Cal K({k,l})$.

If $y\in\Sigma({k,l})$ then  the limit
 of $x^{(j)}$ in     $\Cal K({k,l})$ is defined to be $y$.

{\it The second step.} 
Let $y\notin \Sigma(I_{k,l})$. Then $y$ is a element of some open 
simplex $\widetilde\Delta(I_{k,l}\ |\goth t)$, i.e $y$  has the form
$$
\{ 1=y_k=\dots=y_{\gamma_1}>y_{\gamma_1+1}=\dots=y_{\gamma_2}>
\dots>y_{\gamma_{s-1}+1}=\dots=\gamma_l=0         \}
$$
In this case  the sufficient and nessesary condition of
convergence of the sequence $x^{(j)}$ in $\Cal K({k,l})$
is the convergence of all  sequences
$$
x_{[\psi]}^{(j)}:=(x_{\gamma_\psi+1}^{(j)},\dots,x_{\gamma_{\psi+1}}^{(j)})
\in\Sigma(I_{\gamma_\psi+1,\gamma_{\psi+1}})
$$
in the Karpelevich velocity polyhedra $\Cal K({\gamma_\psi+1,\gamma_{\psi+1}}  )$
(this convergence is defined by the inductive assumption)

          This concludes the definition.

\demo{Example 3.8}
Let $k=1,l=8$.
$$
x_1^{(j)}=2j^3 \ \ \ \ \ \ x_2^{(j)}=j^3  \ \
$$
$$
x_3^{(j)}=j^2+j+2 \ \ \ \ \ x_4^{(j)}=j^2+j+1 \ \ \ \ x_5^{(j)}=j^2+j \ \
$$
$$
x_6^{(j)}=2j \ \ \ \ \ x_7^{(j)}=j \ \ \ \ x_8^{(j)}=0
$$
Then the associated tree-partition has the form

$$
(1\ \ 2 \ \ 3 \ \ 4 \ \ 5 \ \ 6 \ \ 7 \ \ 8) 
$$
$$
( \ 1 \ ) \ \ ( \ 2 \ ) \ \ (3 \ \ 4 \ \ 5 \ \ 6 \ \ 7 \ \ 8)
 $$
$$
\ \ \ \ \ \ \ \ \ \ \ \ (3 \ \ 4 \ \ 5 \ ) \ \ ( 6 \ \ 7 \ \ 8)
$$
$$
\ \ \ \ \ \ \ \ \ \ \ \ \ \ \ \ \ \ \ \ \ (\ 6 \  )  \ \ (\ 7 \ ) \ (\ 8\ )
$$
 The limit of $x^{(j)}$ in $\overline\Sigma(I_{1,8})$ is the 
collection
$$
\{ 1>1/2>0=0=0=0=0=0\}\in \Delta(I_{1,8})
\tag 3.8 
$$
The sequence $x^{(j)}$ induces the sequence
 $$y^{(j)}=(x_3^{(j)},...,x_8^{(j)})\in
\Sigma(I_{3,8})$$
 The limit of $y^{(j)}$  in $\overline\Sigma(I_{3,8})$ is the collection
$$\{1\ge 1\ge 1\ge 0\ge 0\ge 0\} 
\in\Delta(I_{3,8})\tag 3.9$$
Now we obtain the sequences
$$ z^{(j)}=(x^{(j)}_3,x^{(j)}_4,x^{(j)}_5)\in\Sigma(I_{3,5})
$$
$$ u^{(j)}=(x^{(j)}_6,x^{(j)}_7,x^{(j)}_8)\in\Sigma(I_{6,8})
$$
We have
$$z^{(j)}=(j^2+j+2, j^2+j+1 ,j^2+j  )=(2,1,0)
$$
(recall that the collection $z^{(j)}$ is defined up to additive constant) and $\lim z^{(j)}$ is the collection
$$
\{2>1>0\}\in \Sigma(I_{3,5}) \tag 3.10
$$
At last
$$
u^{(j)}=(2j,j,0)
$$
and the limit of $u^{(j)}$ in $\overline\Sigma(I_{6,8})$
is the point
$$
\{1>1/2>0\}\in\widetilde\Delta(I_{6,8})\tag 3.11
$$
The limit of the sequence $x^{(j)}$ is the collection of
 collections (3.8 )-(3.11 ). 
\enddemo
\enddemo

\demo{Topology on the boundary of $I_{k,l}$}
This topology satisfies the following property: the closure of 
$F(\goth a)$ consists of all faces $F(\goth b)$ such that $\goth b 
<\goth a$.

  We assume  the topology is defined for all polyhedra 
  $\Cal K(\alpha,\beta)$ such that $\beta-\alpha<l-k$.

We define the convergence of a sequence 
$$
Z^{(j)}\in F(\goth a)
$$
in two steps.

{\it The first step } Let
$$
h^{(j)}=\{ 1=h^{(j)}_k\ge h^{(j)}_{k+1}\ge\dots\ge h^{(j)}_l=0\}
$$
be the component of $Z^{(j)}$ associated to multiplier 
$\widetilde\Delta(\goth a \ , I_{k,l})$ in the product (3.7 ).
Then the convergence of $h^{(j)}$ in $\Delta(\goth a \ , I_{k,l})$
is a nessesary condition for its convergence in $\Cal K(k,l)$.
We denote  the limit of $h^{(j)}$ in   $\Delta(\goth a \ , I_{k,l})$ 
by $u$.
   
{\it Second step}.  Let us consider the partition of $I_{k,l}$ associated to $\goth a$ :
$$
I_{k,l}=I_{k,\gamma_1}\cup
I_{\gamma_1+1,\gamma_2}\cup\dots\cup I_{\gamma_{s-1}+1,l}
$$
Then the collection $u$ has the form
$$
u=\{ 1=u_k=\dots=u_{\gamma_1}\ge
u_{\gamma_1+1}=\dots=u_{\gamma_2}\ge\dots\}
$$
 Let us consider $\tau_1,\tau_2,\dots$ such that
$$\{1=u_{k}=\dots=u_{{\tau_1}}>u_{\tau_1+1}=\dots=u_{\tau_2} >\dots\}
$$
The set $\{\tau_1,\tau_2,\dots\}$ is a subset in the set 
$\{\gamma_1,\gamma_2,\dots\}$) and hence each segment 
$I_{\tau_\alpha+1},\tau_{\alpha+1}$ is the union of the segments 
$I_{\gamma_m+1,\gamma_{m+1}}$.

   Let us consider on each set
$$\{\tau_\alpha+1,\tau_\alpha+2,\dots,\tau_{\alpha+1}\}
$$
the tree-partition $\goth b_\alpha$ induced by the tree-partition 
$\goth a$. The sequence $Z^{(j)}$ induces the sequence 
$Z^{(j)}_{[\alpha]}$
in each face 
$F(\goth b_\alpha )\subset\Cal K(\tau_{\alpha}+1,\tau_{\alpha+1})$ .

The nessessary and sufficient condition of the convergence of $Z^{(j)}$ is
the convergence of each sequence $Z^{(j)}_{[\alpha]}$ in Karpelevich
polyhedron   $\Cal K(\tau_{\alpha}+1,\tau_{\alpha+1})$.

\enddemo

\demo{3.3. The compactification of symmetric space by Karpelevich 
velocities}
Let us consider the boundary
$$
\partial \Cal K ({1,n}):= \Cal K ({1,n}) \backslash \Sigma(I_{1,n})
$$
of the polyhedron $ \Cal K ({1,n})$.

We define the compactification
$$
(SL(n,\R)/SO(n))\cup (\partial \Cal K ({1,n}))
$$
of the symmetric space $SL(n,\R)/SO(n)$. Let 
$x^{(j)}\in SL(n,\R)/SO(n)$ be a sequence and $y\in\partial \Cal K (I_{1,n})$. The $x^{(j)}\rightarrow y$ if

 1. distance  $d(x^{(j)},0)\rightarrow \infty$

2. $\Lambda(x^{(j)})\rightarrow y $ in the topology of $\Cal K (I_{1,n})$
(where $\Lambda(\cdot) $ is defined by the formula (3.1))
\enddemo\enddemo

\head{4.Tits building on matrix sky}\endhead

We recall that geodesics in the spase          $SL(n,\R)/SO(n)$ have the 
form
$$
\gamma (s)= A
\pmatrix \exp (\lambda_1 s) & \ & \ & \ \\
\ & \exp (\lambda_2 s) & \ & \\
\ & \ & \ddots \ & \\
\ & \ & \ & \ \exp(\lambda_n s)\endpmatrix A^t \tag 4.1
$$
  where $$A\in SO(n),\lambda_1\ge\lambda_2\ge\dots\ge\lambda_n \tag 4.2$$
The term {\it geodesic} below
 means the oriented geodesics without 
fixed parametrization.

\demo{4.1. Matrix sky (visibility boundary)}
Let us consider a Riemann noncompact symmetric space $G/K$.
 Fix a point $x_0 \in G/K$ (in our case $G/K=SL(n,\R)/SO(n)$
it is natural to assume $x_0=E$). Let $T_{x_0}$ be the
tangent space in the point $x_0$ (in our case 
$G/K=SL(n,\R)/SO(n)$ the tangent space is the space
of symmetric matrices defined up to addition of a scalar
matrix, i.e. $A\simeq A+\lambda E$). Let $S$ be the space 
of rays in $T_{x_0}$ with origins in 
zero (i.e. $S=(T_{x_0}\backslash 0)/ \R_+^\ast$ 
where $\R_+^\ast$ is the multiplicative group of positive
real numbers).
  Let  $v\in S$, let $\widetilde v \in T_{x_0}$ be 
a tangent vector on the ray $v$. Let
$$
\gamma_v = \gamma_v(t) $$
be the geodesic such that
$$
\gamma_v(0)=x_0 \ \ \ \gamma_v'(0)=\widetilde v
$$
 We don't interested by the parametrization of the geodesic
$\gamma(s)$ but its direction is essential for us.

 Let  $Sk$ be another copy of the sphere $S$. Points of the sphere $Sk$ 
we consider as infinitely far points of $G/K$. We will call
the sphere $Sk$ by {\it the matrix sky} or by {\it the visibility 
boundary}. Let us
 describe the topology on the space
$$\overline{(G/K)}^{vis}:=G/K\cup Sk$$

We equip the spaces $G/K $ and $Sk $ with the usual 
topology. Let $y_j$ be a sequence in $G/K$. Let $v\in Sk$. Let
$\gamma^{(j)}$ be the geodesic joining points $x_0$
and $y_j$. Let us consider the vectors $v_j\in S$ such that
$$
\gamma^{(j)}=\gamma_{v_j}
$$
The convergence of  the sequence $y_j\in G/K$ to a point $v\in Sk$ 
is defined by the conditions

  1.$ \rho(x_0,y)\rightarrow \infty$

  2.$v_j\rightarrow v$ in the natural topology of the sphere
$Sk$
\enddemo

\demo{4.2. The projection of the matrix sky to the velocity 
simplex}
 Let $G/K=SL(n,\R )/SO(n)$. Let us consider a geodesic $\gamma $
with the origin in $x_0=E$. Then $\gamma $ has the form
$$
\gamma (s)= A
\pmatrix \exp (\lambda_1 s) & \ & \ & \ \\
\ & \exp (\lambda_2 s) & \ & \\
\ & \ & \ddots \ & \\
\ & \ & \ & \ \exp(\lambda_n s)\endpmatrix A^t \tag 4.3
$$
where $A\in SO(n)$ and
$$
\lambda_1\ge \lambda_2\ge ...\ge\lambda_n=0$$
Let $\Delta=\Delta_n$ be the simplex
$$
1\ge\mu_2\ge ...\ge\mu_{n-1}\ge 0
$$
(see 3.1) . We associate to each geodesic $\gamma(s)$ the
point 
$$
D(\gamma) : 1\ge {\lambda_2 \over \lambda_1} \ge{
\lambda_3\over \lambda_1}\ge ...\ge {\lambda_{n-1} \over \lambda_1}\ge 0
$$
of the simplex $\Delta$ .

Obviously $D(\gamma)$ is the limit of the geodesics $\gamma$ in 
the simplest velocity compactification of $SL(n,\R)/SO(n)$.
We say that $D(\gamma)\in\Delta $ is the {\it velocity of geodesic}
$\gamma$

\enddemo

\demo{4.3. The projection of the matrix sky to the space of flags}
Let $\Cal F$ be the set of all flags
$$
0=V_0\subset V_1\subset V_2 \subset ... \subset V_s=\R^n$$
in $\R^n \ (s=0,1,2,\dots, n)$, see section 7.  Denote by 
$\Cal F_{complete}$ the space of complete flags (i.e $i=n$)

   Let us consider the geodesic $\gamma (s)$ given by the expression
(4.1). Let the collection $\lambda_1, \lambda_2,...,\lambda_n$
 has the form
$$
\lambda_1=\lambda_2=\dots =\lambda_{s_1}>\lambda_{s_1+1}=
 \lambda_{s_1+2}=\dots =\lambda_{s_2}>\dots      \tag 4.4
$$
Let $T_\alpha$ be the subspase in $\R^n$ which consists
of vectors
$$
(x_1,\dots,x_{s_\alpha},0,0,\dots )
$$
Let $V_\alpha=AT_\alpha$ (see (4.3)). We denote by $F(\gamma)$
the flag
$$
V_1\subset V_2 \subset V_3 \subset\dots \tag 4.5
$$
We obtain the map $F:Sky\rightarrow \Cal F$. 
It is easy to see that the geodesic $\gamma$ is determined by the pair
$$
(D(\gamma); F(\gamma))\in\Delta \times \Cal F
$$
A pair (velocity (4.4), flag (4.5)) is not arbitrary. It has to 
satisfy the condition $\dim V_j=s_j$.
\enddemo
\demo {4.4. Limits of geodesics on the matrix sky}
Let us consider arbitrary geodesic $\gamma(s)$given by the formula 
(4.1)-(4.2). Let us consider the geodesics $\kappa_s(t)$ joining the 
points $0$ and $\gamma(s)$. We want to calculate $\lim_{s\rightarrow\infty } 
\kappa(s)$. 

  For this purpose let us represent  the matrix $A\in GL(n,\R)$ in 
  the form $A=UB$ where $U\in O(n)$ and $B$ is uppertriangle matrix. It 
  is
easy to prove that the limit of the family of geodesics $ \gamma_s$
is the geodesics $\sigma(t)$ given by the formula

$$
\sigma(t)= U\pmatrix \exp(\lambda_1 t)&&\\
    &\ddots &     \\
&&\exp(\lambda_n t)\endpmatrix U^{-1}
$$\
This remark has several simple corollaries          .

\demo  {  Construction of the matrix sky doesn't depend on the point 
$x_0$}

Indeed let us consider two points $x_0$ and $x_1$ and denote the associated 
matrix skies by $Sk(x_0), Sk(x_1)$          . Let us consider a geodesic 
$\gamma(s)$ with the origin in $x_1$. Then $\gamma(s)$ has limit on
$Sk(x_0)$ . Hence we obtain the canonical map
 $\psi_{10}:Sk(x_1)\rightarrow Sk(x_0)$   . We also have canonical 
 map $\psi_{01}: Sk(x_0)\rightarrow Sk(x_1)$.
It is easy to show that 
$\psi_{01}\circ\psi_{10}=id,\psi_{10}\circ\psi_{01}=id$  and we
obtain the canonical bijection $Sk(x_0)\leftrightarrow  Sk(x_1)$.
               \enddemo
\demo { In particular               for each point $x\in G/K$ and each 
point $y\in Sk$ there exists the unique geodesic joining $x$ and 
$y$    }
   \enddemo
\demo{ The group $G/K$ act by the natural way on the space
$(G/K)^{vis}$}
     
 Indeed the group $G$ acts on the space of geodesics.
     \enddemo
For each $g\in G$ and each $\gamma\in Sk$
$$
D(g\cdot \gamma)=D(\gamma) \ \ \ \ F(g\cdot\gamma)=g\cdot F(\gamma)
$$                       
 \enddemo

\demo{4.5.
Simplicial structure on the matrix sky}
Let us consider a complete flag    $L\in\Cal F_{complete}$
$$
L: 0\subset W_1 \subset W_2 \subset\dots \subset W_{n-1}\subset \R^n;\ \ \
\dim W_j=j
$$
Let us consider the embedding
$$
\sigma_L:\Delta \rightarrow Sk
$$
defined by the conditions

1. $G\circ \sigma_L $ is the identity map $\Delta\rightarrow
\Delta$

2. The image of the map $F\circ \sigma_L:\Delta\rightarrow \Cal F $ 
consists of subflags of the flag $L$.

 Now we will give a explicit construction of the map 
$\sigma_L$. Without loss of generality we can consider
the flag
$$
\R^1\subset\R^2\subset \dots \subset \R^{n-1}
$$
in $\R^n$, the subsubspace $\R^j$ consists of vectors
$(x_1,\dots ,x_j , 0 , \dots , 0)$. Let
$$
1\ge\mu_2\ge ...\ge\mu_{n-1}\ge 0
$$
be a point of $\Delta$. Then the associated geodesic (we remind that 
geodesic  
is identified with the point
of $Sk$) has the form
$$\gamma(s)=
\pmatrix \exp(s) &  &  & &\\
& \exp(\mu_2s)& & &\\
& & \ddots & &\\
&  & & \exp(\mu_{n-1}s) &\\
& & & & 1 \endpmatrix
$$
 Hence we obtain the tiling of the sphere $Sk$ by 
the simplices $\sigma_L(\Delta)$. These simplices are enumerated by 
the points      $L$
of the spase of {\it complete}
flags. It is easy to show that this tiling satisfies the conditions

  a) Let $g\in SL(n,\R)$. Then
$$
\sigma_{gL}(\Delta)=g\cdot \sigma_{L}(\Delta)$$

  b) If $L\not=L'$ then the interiors of simplices
$\sigma_{L}(\Delta)$ and $\sigma_{L'}(\Delta)$ doesn't intersect

  c) Let 
$$
L:V_1\subset V_2 \subset \dots V_{n-1}
$$
$$
L':V_1'\subset V_2' \subset \dots V_{n-1}'
$$
be complete flags. If $V_j\not=V_j'$ for all $j$ then
$\sigma_{L}(\Delta)\cap \sigma_{L'}(\Delta)=\emptyset$ .
In the opposite case
 the intersection
$$\Lambda=\sigma_{L}(\Delta)\cap \sigma_{L'}(\Delta)$$
is a joint face of simplices $\sigma_{L}(\Delta)$ and $\sigma_{L'}(\Delta)$. 
Let us describe $\Lambda$. Let $\alpha_1,\dots,\alpha_s$ be all 
indices $j$ such that
$V_j=V_j'$ (i.e $V_{\alpha_i}=V_{\alpha_i}'$ and $V_j\not=
V_j'$ for all $j\not=\alpha_i$).  Let us consider the face
$$
1=\lambda_1=\dots=\lambda_{\alpha_1}>\lambda_{\alpha_1+1}=
\lambda_{\alpha_1+2}=
\dots=\lambda_{\alpha_2}>\dots$$
of the simplex $\Delta$. Then
$$\Lambda=\Sigma_L(N)=\Sigma_{L'}(N)$$
 Now we  obtain on the sphere $Sk$ the structure of a Tits 
building (see [30])
\enddemo

  \demo{4.6. Tits metric on the matrix sky}
Let us consider points $y_1, y_2 \in \sigma_{L}(\Delta)$. We define 
the distance $d(y_1, y_2 )$ as the angle between geodesics
$x_0y_1$ and $x_0y_2$.
  Let $z,u \in Sk$. Let us consider a chain
 $$
z=z_1,z_2,\dots, z_\beta=u \ \ \ \ (z_j\in Sk)
$$
such that for all $j $ points $  z_J,z_{j+1}$ 
belongs to one element of our tiling.

 Let us define the {\it Tits metric} $D(\cdot,\cdot)$ on $Sk$by the formula
 $$
D(z,u)=\inf(\sum_j d(z_j),d(z_{j+1})  )
$$
(we consider the infimum   by the all chains $z_1,\dots,z_\beta$).

 \demo {Remark 4.1  }
The topology on the $Sk$ defined by the Tits metric is not
equivalent to the usual topology of the sphere.
\enddemo

\demo {Example 4.2 }
Let $n=3$,$G/K=SL(3,\R)/SO(3)$. Then $Sk$ is the 4-dimensional 
sphere $S^4$, $\dim\Delta=1$,i.e the simplices
$\sigma_{L}(\Delta)$ are segments. We will describe the siplicial 
structure on $Sk=S^4$. Let $P$ be the spase of all 1-dimensional 
linear subspaces in $\R^3$ and $Q$ be the space of 2-dimensional 
subspaces in $\R^3$ (evidently
$P\simeq Q$ are the projective planes). We 
want to construct some graph $\Gamma$. The set of vertices of $\Gamma$
is $P\cup Q$. Let $p\in P,q\in Q, p\subset q$. Then $p$ and $q$ are adjacent
 to the same edge and all edges have such form. Assume that the length of each 
edge is $\pi/3$. Then
graph $\Gamma$ is isometric to the sphere $Sk=S^4$ endowed with the 
Tits metric.

\enddemo
\enddemo

\demo{4.7. Abel subspaces}
Let $A$ be a orthogonal matrix. Let us consider the submanifold 
$R[A]\subset SL(n,\R)/SO(n)$ consisting of matrices of the form

$$
\psi_A(s_1,\dots,s_{n-1})=A
\pmatrix \exp(s_1)& & & &\\
& \exp(s_2) &&&\\
&&\ddots&&\\
&&&\exp(s_{n-1})&\\
&&&&1\endpmatrix A^{-1}
$$
where $s_1,\dots,s_{n-1}\in\R$.

 The map
$$(s_1,\dots,s_{n-1})\mapsto \psi(s_1,\dots,s_{n-1})
$$
is the isometric embedding $\R^{n-1} \rightarrow SL(n,\R)/SO(n)$ (with 
respect to the standard metrics
in $\R^{n-1}$ and in $SL(n,\R)/SO(n)$).

  Let us consider the trace $S[A]$ of the space $R[A]$ on the
surface $Sk$. It is easy to see that $S[A]$ is the union
of $(n-1)!$ simplexes $\sigma_{L}(\Delta)$. These simplices
are separated by the hyperplanes $s_i=s_j$.

 \enddemo 

\head{5. Hybridization: Dynkin-Olshanetsky and Karpelevich boundaries}
\endhead

\demo{5.1 Hybridization}
  Let
$$
i_1:G/K\rightarrow X
$$
$$
i_2:G/K\rightarrow Y
$$
be embeddings of symmetric space $G/K$ to compact metric spaces $X$ and $Y$.
Let the images of $G/K$ in $X$ and $Y$ be dense.
 
 Let us consider the embedding

$$
i_1\times i_2:G/K\rightarrow X\times Y
$$
defined by the formula
$$
h\mapsto (i_1(h),i_2(h))
$$
where $h\in G/K$. Let $Z$ be the closure of the image of $G/K$ in $X\times Y$.
Then $Z$ is the new compactification
 of $G/K$. We say that $Z$ is the {\it hybrid} of
$X$ and $Y$.

  We want to apply this construction in the case then $X$ is a velocity
compactification and $Y$ is Satake-Furstenberg compactification.
\enddemo

\demo{5.2. Dynkin-Olshanetsky boundary}
Let us consider the hybrid $Z$ of the simplest  velocity compactification
(see 3.1) and Satake-Furstenberg compactification of Riemann noncompact
symmetric space. Again let us consider only the case $G/K=SL(n,\R)/O(n)$.

  A point of the space $Z$ is given by the following data

$0^\ast$. $s=1,2,...,n-1$

$1^\ast$. A hinge
$$
\goth P= (P_1,\dots,P_s)$$
such that  $P_j$ are nonnegative  definite
(see section 9     )

$2^\ast$. A point of the simplex $\Delta_s$:
$$
1\ge\mu_2\ge\dots\ge \mu_{s-1}\ge 0
$$

  Let $x^{(j)}\in SL(n,\R)/O(n)$ be a unbounded
 sequence. Let $a_1^{(j)}\ge\dots
\ge a_n^{(j)}$ be the eigenvalues of
 $x^{(j)}$. Let $\lambda_\alpha^{(j)}=\ln a_\alpha^{(j)}$. Then 
the point $\Lambda(x^{(j)}):= (\lambda_1^{(j)},\lambda_2^{(j)},\dots)$ be a
 point
of the the simplicial cone $\Sigma_n$(see 4.1).
 The sequence $x^{(j)}\in SL(n,\R)/O(n)$ converges
in $Z$ if $x^{(j)}$ converges in Furstenberg-Satake compactification and
$\Lambda(x^{(j)})$ converges in the velocity simplex $\overline\Sigma_n=
\Sigma_n\cup \Delta$.

Now we want to explain how to calculate $\lim x^{(j)}$.
  Let $\goth P= (P_1,\dots,P_s)$ be the limit of $x^{(j)}$ in 
Satake-Furstenberg
compactification. Let $\gamma_j=\dim Im (P_j)$. 
Let $(\tau_2,\dots,\tau_{s-1})$
be the limit of $\Lambda(x^{(j)})$ in the simplex $\Delta$. Then the 
collection $\tau_2\ge\tau_3\ge\dots$ has the form
$$
1=\tau_1=\dots=\tau_{\gamma_1}>\tau_{\gamma_1+1}=\dots=\tau_{\gamma_2}>\dots
\tag 5.1
$$
We assume
$$\mu_j:=\tau_{\gamma_{j-1}+1}=\dots=\tau_{\gamma_j}\tag 5.2 $$
and we obtain the data $0^\ast$ - $2^\ast$.
\enddemo

\demo{5.3. 
The projection of the Dynkin-Olshanetsky boundary to the 
matrix sky}
Let we have data $0^\ast$ - $2^\ast$
Let us consider the new data

$1^+$. The flag
$$
Ker(P_1)\supset Ker(P_2)\supset Ker(P_3)\supset\dots
$$
$2^+$. The collection of numbers $\tau_2,\dots\tau_{s-1}$ defined 
by the formula (5.2)

 These data define the point of the matrix sky (see  4.2-4.3  )
\enddemo

\demo{5.4. Limits of geodesics}
Let us consider a geodesics
$$
\gamma(s)=A
\pmatrix \exp(\lambda_1 s) & &&\\
&\exp(\lambda_2 s) &&\\
&&\dots &\\
&&&1\\ \endpmatrix
A^t
$$
where $A\in SL(n,\R)$ and $\lambda_1\ge\lambda_2\ge\dots\ge\lambda_n= 0$.
 Let $\goth P= (P_1,\dots,P_s)$ be the limit of $\gamma(s)$ in the 
space of hinges. The limit
of $\gamma(s)$ in the velocity simplex $\Delta_n$ is 
$\tau_2,\dots,\tau_{n-1}$
where $\tau_j=\lambda_j/\lambda_1$.
 
   Let $\gamma_\alpha= \dim Im(P_\alpha)$. We define numbers
$$
\mu_\alpha:=\tau_{\gamma_{\alpha-1}+1}=\dots=\tau_{\gamma_{\alpha}}
$$
Now we  obtain the data $0^\ast$ - $2^\ast$.
\demo{Remark 5.1} Not all points of Dynkin-Olshanetsky boundary are
limits  of geodesics. A point defined by the data $0^\ast-2^\ast$
is the limit of a geodesics if and only if $1>\mu_2>\dots>\mu_{n-1} >0$
\enddemo
\enddemo

\demo{5.5. Karpelevich compactification}
The Karpelevich compactification is the hybrid of 
 the compactification by
 Karpelevich velocities and Satake-Furstenberg compactification.

  A point of the Karpelevich compactification is givem by the following 
 data

$0^\star$. $s=1,\dots,n-1$

$1^\star$. A hinge
$$\goth P= (P_1,\dots,P_s)$$
such that  are positive definite (see  section 10       )

$2^\star$. A point of the  boundary of the Karpelevich
 velocity polyhedron $\Cal K(1,s)$ (see 3.2)

   The topology on Karpelevich compactification is defined by the obvious 
   way
.

   The natural projection $\partial\Cal K(1,s) \rightarrow 
\Delta(I_{1,s})$ defines the
projection of Karpelevich boundary to Dynkin-Olshanetsky boundary

\enddemo

\head{6. Space of geodesics and sea urchins}\endhead

\demo {6.1. Spase of geodesics}
Let us consider a Riemann noncompact symmetric space $G/K=SL(n,\R)/SO(n)$.
 Denote by $\goth G$
the spase of all oriented geodesics in $G/K$.

The question about topologies on $\goth G$ is delicate. I'll describe 
the topology which seems to me the most natural.

Let us consider a collection of integers $A=(\alpha_0,\dots,\alpha_\sigma)$
such that
$$
1=\alpha_0\le\alpha_1\le\dots\le\alpha_\sigma=n
$$
Let us denote by
$\Delta(A)$              the open simplex
$$
1=\lambda_1=\dots=\lambda_{\alpha_1}>\lambda_{\alpha_1+1}=\dots=
\lambda_{\alpha_2}>\dots>\lambda_{\alpha_{\sigma-1}+1}=\dots=\lambda_n=0
$$
Simplices $\Delta(A)$ don't intersects and $\cup_A\Delta(A)$
coincides with the simplex $\Delta_n$.

Let us consider a geodesic $\gamma\in\goth G$. Its velocity is a point of one of the 
simplices $\Delta(A)$. The space of all geodesics with a given 
velocity $\Lambda\in\Delta(A)$ is   a $SL(n,\R)$-homogeneous space.
The stabilizer $G(A)$ of the geodesic $\gamma$ (up to conjugacy)
 depends only of the the collection $A$ (it doesn't depend of $\Lambda$
and the geodesic itself):
$$
G(A)=\R_+^\ast \times \prod O(\alpha_{j+1}-\alpha_j).
$$
We denote by $\goth G(A)$ the space of all geodesics which velocities 
are elements of $\Delta(A)$. Then
$$
\goth G(A)\simeq \Delta(A)\times (SL(n,\R)/G(A))
$$
We equip this space with the usual  topology of the direct
product. We equip the space
$$
\goth G=\bigcup_A   \goth G(A)\simeq 
\bigcup_A \Delta(A)\times (SL(n,\R)/G(A))
     $$
with  the topology of disjoint union.
 \demo{Remark 6.1}Hence the space of geodesics is disconnected set.
It is not strange.        Let $A_0=\{0,1,2,\dots,n\}$
Let us consider the set of limits of the geodesics 
$\gamma\in \goth G(A_0)$ on
matrix sky. Then this set is open and dence. The set of limits
 of $\gamma\in \goth G(A)_0$
in Satake boundary is compact. Hence it is natural to think that
$\goth G $ is disconnected.
\enddemo
\enddemo

\demo{6.2. Spase of geodesics as boundary of symmetric space}
 We define the natural topology on
 the space 
$$\goth R=G/K \cup \goth G$$
We equip the space $G/K $ with the
natural topology. The space $\goth G$ is 
 equipped with the 
topology mentioned above and the spase $\goth G $ 
is closed in $\goth R$. Fix a point
 $b_0\in G/K$. Let $x_j\in G/K$ be a unbounded sequence. 
The seguence $x_j$ converges in $\goth R$
if it  satisfies the 
following conditions 

  1. Sequence of geodesics $b_0 x_j$ converges. Denote 
by $y$ its limit on the matrix sky.

  2. There exists a limit $z$ of the sequence of geodesics $yx_j$.

The limit of the sequence $x_j$
is defined to be the geodesic $z$.

\demo{Remark 6.2} In our case the dimension of the boundary
$$
\dim \goth G = 2 \dim G/K -2
$$
is greater than $\dim G/K$ (even in the case then $G/K=SL(2,\R)/SO(2)$ is 
Lobachevskii plane)
\enddemo

\demo {Remark 6.3}
The spase $\goth R $ is not compact(since $\goth G$ is not compact)
\enddemo
\enddemo

\demo{6.3. Sea urchins}
Recall that each geodesic
 $\gamma\in\goth G$ has a velocity 
$\{\mu_2,\mu_3,\dots\}$ which is a point
of the simplex $\Delta$ (see 3.1). We denote by $\goth G^{rat}$ 
the space of geodesics having rational velocities (i.e. 
$\mu_j$ are rational). Let us consider the 
set ({\it sea urchin}                              )
$$ 
\goth R^{rat}:=G/K \cup \goth G^{rat} \subset \goth R
$$
We don't interested by the topology on
sea urchin (it is seems natural to consider the discrete topology on the 
set of velocities, the usual topology
on the space of geodesics with a given velocity and the natural (see 
6.2)
convergence of sequences in $G/K$ to geodesics)
\enddemo

\demo{6.4. Projective universality}
Let $\rho_j$ be a finite family of linear irreducible representations of
 the group $G$ in the spaces $V_j$. We assume that for each $j$ there
 exists a $K$-fixed nonzero vector
$v_j \in V_j$.
Let us consider the direct sum $\rho=\oplus \rho_j$
   of representations $\rho_j$ and the vector $w=\oplus v_j\in\oplus V_j$. 
Let us consider the projective space $\P (\oplus V_j)$
.
 Let $\Cal O\simeq G/K$ be the $G$-orbit of the vector 
 $w\in \P(\Cal O )$.
 Let $\overline {\Cal O}$ be the closure of $\Cal O$ in
$\P (\oplus V_j)$.

The $G$-spaces $\overline {\Cal O}$ are called 
{\it projective compactifications } of $G/K$
 
   We will construct the map
$$\pi:\goth R=G/K \cup \goth G \rightarrow \overline {\Cal O}
$$
The map $G/K\rightarrow \Cal O$
 is obvious. Let us consider a geodesic $\gamma (s)\in\goth G$.
 It is easy to prove that there exists
 $\lim_{s\rightarrow\infty}\rho(\gamma(s))$
in $\P (\oplus V_j)$. By definition $\pi (\gamma) $ is this limit.
 
  \proclaim{Proposition 6.4}a) The map $\pi:\goth R \rightarrow 
\overline {\Cal O}$
is surjective.

  b) Moreover the $\pi$-image of sea urchin $\goth R^{rat}$
is the whole $\overline {\Cal O}$.
\endproclaim

\enddemo

\head{7.Bibliographical remarks}
\endhead

\demo{Remarks to section 1-2}
The Satake-Furstenberg boundary is a version of Study-Semple-Satake-
Furstenberg-De Concini-Procesi-Oshima boundary (see [1-7]) 
of symmetric space
$G/H$ where $G$ is a semisimple group and $H$ is a symmetric subgroup
(i.e. subgroup $H$ is set of fixed points for some involution on the
group $G$ ). The usual definition is the following. Let us consider
a finite-dimensional irreducible representation of $G$ having
a $H$-fixed vector $v$(the representation $\rho$ have to satisfy some
nondegeneracy conditions). Then our compactification is the closure 
of the orbit $G\cdot v$ in the projective space. The coincidence of 
our construction with classical is not obvious, for construction
of projective embedding of space $Hinge(n)$ see [8,10].

  Hinges were defined in [8], see also [10]. For construction of 
  separated quotient space through Hausdorff metric see [9].
For construction of separated quotient space it is also possible
to use closure in Chow scheme , see [11-12].

  Our space $Hinge(n)$ is one of the real form of Semple complete collineation 
  variety. The Satake-Furstenberg compactification of $SL(n,\R)$
is one of the real forms of Study-Semple complete quadrics.

Data $1^\star-2^\star$  were introduced in [3].
 \enddemo
\demo {Remarks to section 3}
I haven't seen this construction in literature.  The analogy
of the collection$\{\ln \lambda_J\}$ for arbitrary symmetric space
is so-called complex (or compound) distance (see for instance [10])

 Kaprelevich velocity
polyhedron is the closure of a Weyl chamber in Kaprelevich 
compactification.
       \enddemo

\demo{Remarks to section 4} See [15,32].

The most of constructions described in this paper are very exotic 
from the point of view of the official differential geometry. The 
visibility boundary is exeption. It is more or less general 
differential-geometric object, see [13-15].

   Tits metric on the infinitely distant sphere(see also [16] for 
   boundaries of Bruhat-Tits buildings) also is more or less
general construction (see [15]). Neverless nice tiling of the sphere
also seems exeptional phenomena.
\enddemo
\demo {Remarks to section 5}.
Karpelevich boundary was constructed in [17] in terms of geometry of 
geodesics. Dynkin-Olshanetsky boundary (see [18-20]) is Martin boundary 
(see [22-25]) for the diffusion on symmetric spacesDiscussion
of these boundaries see also [21].
        \enddemo
\demo {Remarks to section 6}
I havn' seen sea urchin construction in literature. See [26-27]
for universal projective compactification of symmetric space
(see urchin is not compact).

                \enddemo

\enddocument